\documentclass[twoside]{amsart}
\usepackage{amssymb,amsmath}
\usepackage[all]{xy}
\usepackage{fancyhdr}

\def\bc{\begin{center}}
\def\ec{\end{center}}
\def\no{\noindent}
\def\lra{\Leftrightarrow}

\def\Ra{\Rightarrow}
\def\ra{\rightarrow}
\def\ol{\overline}

\def\suse{\subseteq}

\def\im{{\rm Im}}
\def\ker{{\rm ker}}

\def\cok{{\rm cok}}
\def\Ext{{\rm Ext}}
\def\Tor{{\rm Tor}}

\def\pd{{\rm pd}}
\def\id{{\rm id}}
\def\fd{{\rm fd}}
\def\cd{{\rm cd}}

\def\dim{{\rm dim}}

\def\FPD{{\rm FPD}}
\def\FFD{{\rm FFD}}

\def\Hom{{\rm Hom}}

\def\lst{\leqslant}
\def\gst{\geqslant}

\def\lid{{L_n{\rm id}}}
\def\ld{{L_n\dim}}
\def\gl{{\rm gl.dim}}
\def\wgl{w.{\rm gl.dim}}

\newtheorem{thm}{Theorem}[section]
\newtheorem{defn}[thm]{Definition}
\newtheorem{cor}[thm]{Corollary}
\newtheorem{lem}[thm]{Lemma}
\newtheorem{exa}[thm]{Example}
\newtheorem{pro}[thm]{Proposition}

\begin{document}

\bc{\bf\Large On $L_{n}$-Injective Modules and $L_{n}$-Injective Dimensions}\ec

\vskip3mm
\vskip3mm \bc{Tao Xiong${}^1$,\ \ Fanggui Wang${}^2$,\ \ Lei
Qiao${}^3$,\ \ Shiqi Xing${}^4$,\ \ Qing Li${}^5$}\ec

\bc{1,2,3,4:\ \ College of Mathematics and software science, Sichuan
Normal University, Chengdu, 610068 P. R. China;\ \ \newline 5:\ \
College of Computer Science and Technology, Southwest University for
Nationalities, Chengdu, 610041 P. R. China}\ec

\bc{1:\ \ xiangqian2004@163.com;\ \ 2:\ \ wangfg2004@163.com;\ \ 3:\
\ qiaolei5@yeah.net;\ \ \newline 4:\ \ sqxing@yeah.net;\ \ 5:\ \
lqop80@163.com}\ec

%\thanks{This work was completed with the support of}

\vskip5mm
\begin{abstract}
Let $R$ be a ring, and $n$ a fixed nonnegative integer. An $R$-module $W$ is called $L_{n}$-injective if $\Ext_{R}^{1}(M,W)=0$ for
any $R$-module $M$ with flat dimension at most $n$. In this paper, we prove first that ($\mathcal{F}_{n},\mathcal{L}_{n}$) is a complete hereditary
cotorsion theory, where $\mathcal{F}_n$ (resp. $\mathcal{L}_n$) denotes the class of all $R$-modules
with flat dimension at most $n$ (resp. $L_{n}$-injective $R$-modules). Then we introduce the $L_{n}$-injective dimension of
a module and $L_n$-global dimension of a ring.
Finally, over rings with weak global dimension $\leq n$, perfect rings, and $L_n$-hereditary rings,  more properties and
applications of $L_{n}$-injective modules, $L_{n}$-injective dimensions of
modules and $\mathcal{L}_{n}$-global dimensions of rings are given.
\end{abstract}

\no{\bf Keywords:}\ \ flat dimension; $L_{n}$-injective module; $L_{n}$-injective dimension, $L_{n}$-global dimension, $L_n$-hereditary ring.

\no{\bf MSC (2010):}\ \ 16E05, 16E10

\vskip8mm
\section{Introduction}

\vskip3mm
Throughout this paper, $R$ is an associative ring with identity and all modules are left $R$-modules unless otherwise stated.
We also use ${}_R{\frak M}$ to denote the category of left $R$-modules, $\wgl(R)$ (resp. $\gl(R)$) to denote the weak global (resp. global)
dimension of $R$, $\mathcal{F}_n$ to denote the class of all $R$-modules with flat dimension at most $n$.
For an $R$-module $M$, $\pd_{R}M$ (resp. $\fd_{R}M$) stands for the projective (resp. flat) dimension of $M$,
$\id_{R}M$ (resp. $\cd_{R}M$) stands for the injective (resp. cotorsion) dimension of $M$.

Cotorsion modules have being received a lot of attension in many articles, see \cite{BBE01, En84, Xu96}. As in \cite{En84}, an $R$-module
$C$ is called cotorsion if $\Ext_R^1(F,C)=0$ for all flat module $F$. A celebrated result that was proved by Bican et al. in  \cite{BBE01}
is the Flat Cover Conjecture (FCC): Over any ring, every module has a flat cover and hence every module has a cotorsion envelope.

On the further development on the idea of cotorsion module notion, in \cite{Le06}, the weak-injective modules have been studied by Lee.
Recall from \cite{Le06} that An $R$-module
$W$ is called weak-injective if $\Ext_R^1(M,W)=0$ for all modules $M$ with $\fd_{R}M \leq 1$ and from  \cite{BS03} that
a domain $R$ is called almost perfect (APD shortly) if all its proper
homomorphic image are perfect. It was proved in \cite[Corollary 6.4.8]{FL09} that a domain $R$ is an APD if and only if every
module of flat dimension $\leq 1$ has projective dimension $\leq 1$; if and only if every divisible module is weak-injective; if and only if every epic image
of a weak-injective module is weak-injective.

In 2012 the notion of $n$-cotorsion modules was introduced in \cite{EH12} by Enochs and Huang. An $R$-module $N$ is called $n$-cotorsion in \cite{EH12}
if $\Ext_{R}^{1}(M,N)=0$ for all $R$-module $M$ with $\fd_{R}M\leq n$. But the name of the $n$-cotorsion module has been used by Mao and Ding. In \cite{MD06}
the $n$-cotorsion module $N$ means $\Ext_{R}^{n+1}(F,N)=0$ for any flat $R$-module $F$. The two notions of $n$-cotorsion modules are
not coincident, see Example \ref{e202}.

In this paper, the $n$-cotorsion modules which are defined in \cite{EH12}, following Lee's idea, are said to be $L_n$-injective modules.
Thus weak-injective modules are exactly $L_{1}$-injective modules.
Denote ${\mathcal F}_n$ and ${\mathcal L}_n$ the classes of modules of flat dimension $\leq n$
and of $L_n$-injective modules, respectively.

In Section 2, we prove in Corollary \ref{c207} that the $n$-th cosyzygy of a cotorsion $R$-module $M$ is $L_{n}$-injective, and in
Theorem \ref{t208} that an $R$-module
$M$ has flat dimension $\leq n$ if and only if $\Ext_{R}^{i}(M,N) = 0$ for any $L_{n}$-injective module $N$ and any $i\geq 1$,

Given two classes $\mathcal{A}$ and $\mathcal{B}$ of $R$-modules, set
$\mathcal{A}^{\bot}=\{B\,|\,\mbox{$\Ext_{R}^{1}(A,B)=0$ for all $A\in\mathcal{A}$}\}$
and $^{\bot}\mathcal{B}=\{A\,|\,\mbox{$\Ext_{R}^{1}(A,B)=0$ for all $B\in \mathcal{B}$}\}$, which are called the
right orthogonal class of $\mathcal{A}$ and the left orthogonal class of $\mathcal{B}$, respectively. A pair $(\mathcal{A,B})$ of $R$-modules is
called a cotorsion theory (or cotorsion pair) \cite{EJ00} if $\mathcal{A}^{\bot}=\mathcal{B}$ and $\mathcal{A}=^{\bot}\mathcal{B}$.
A cotorsion theory $(\mathcal{A,B})$ is said to be complete \cite{Tr00} if every $R$-module has a special $\mathcal{A}$-precover.
It is shown in \cite{BBE01} that the pair $(\mathcal{F},\mathcal{C})$ is a complete cotorsion theory, where $\mathcal{F}$ and $\mathcal{C}$
are the classes of flat $R$-modules and cotorsion modules, respectively. In 2006, it is shown in \cite{Le06} that
the pair $(\mathcal{F}_{1},\mathcal{W})$ is a cotorsion theory when $R$ is a domains, where $\mathcal{F}_1$ and $\mathcal{W}$
denote the classes of $R$-modules with flat dimension at most $1$ and weak-injective modules, respectively.
In Section 3, for the further examination, in Theorem \ref{t306}, we prove that the pair $(\mathcal{F}_{n},\mathcal{L}_{n})$ is a complete cotorsion theory,
where $\mathcal{F}_n$ and $\mathcal{L}_n$ are the classes of $R$-modules with flat dimension at most $n$ and $L_{n}$-injective modules, respectively.

In Sections 4 and 5, we are going to introduce the
$\mathcal{L}_{n}$-injective dimension of modules and the $\mathcal{L}_{n}$-global dimension of rings. In Section 6, we start by discussing
when $L_n$-injective modules are injective. It is shown in Theorem \ref{t601} that $L_n$-injective modules are injective if and only if $\wgl(R)\leq n$.
Then we discuss, for $n\geq 1$, when every $R$-module is $L_n$-injective.

 It was shown in \cite{Xu96} that every module is cotorsion if and
only if $R$ is left perfect. In Theorem \ref{t605},  we show that every module is $L_n$-injective if and only if $R$ is left
perfect with $l.{\rm FFD}(R)=0$,
where $l.{\rm FFD}(R)$ is the left finitistic weak dimension of $R$ defined in \cite{Ba60}. We introduce also the notion of $L_n$-hereditary
rings and give a series of consideration on them. Compared the notions of almost perfect rings and $L_1$-hereditary rings, we
point out that all almost perfect rings are $L_1$-hereditary and a domain $R$ is an APD if and only if $R$ is $L_1$-hereditary,
and in Example \ref{e704}, we give an example that some $L_1$-hereditary ring is not perfect.

\vskip5mm
\section{$L_n$-Injective modules}

\vskip3mm
We start this section with the following difinition.

\begin{defn}\label{d20l}{\rm (see \cite[Definition 2.7(1)]{EH12})}
An $R$-module $W$ is called $L_{n}$-injective $(L$ means Lee$)$ if $\Ext_{R}^{1}(M,W)=0$ for all $R$-modules $M\in {\mathcal F}_n$.
\end{defn}

Naturally, $L_0$-injective modules are exactly cotortion modules and $L_1$-injective modules are exactly weak-injective modules.

We denote the class of all $L_{n}$-injective $R$-modules by $\mathcal{L}_n$.

\begin{exa}\label{e202}
It is easy to see that $L_n$-injective modules are $n$-cotorsion modules under the definition of Mao et al. But $n$-cotorsion modules are
not necessarily $L_n$-injective modules. For example, take $n=1$ and let $R=Z$ and $M=R/(2)$.
Then $R$ is a $1$-cotorsion module since $\id_{R}R =1$. It is clear that $\fd_{R}M=1$ and $\Ext_R^1(M,R)\cong R/2R\not=0$ by
{\rm \cite[Theorem 7.17]{Ro79}}. Hence then $R$ is not $L_{1}$-injective.
\end{exa}

\begin{exa}\label{e203}
The following facts are true.

$(1)$\ \ Injective modules are $L_n$-injective modules for all integer $n\geq 0$.

$(2)$\ \ If $m\leq n$, then $L_{n}$-injective $R$-modules are $L_{m}$-injective, and hence $\mathcal{L}_m\supseteq \mathcal{L}_n$.

$(3)$\ \ Let $\{W_i\}$ be a family of $R$-modules. Then $\prod\limits_iW_i$ is $L_{n}$-injective if and only if
each $W_i$ is $L_{n}$-injective modules.

$(4)$\ \ Let $n\geq 1$ and let $R$ be a domain. Suppose $W$ is an $L_n$-injective module. For $a\in R$ and $a\not=0$, since
$\Ext_R^1(R/Ra,W)=0$, we have that $W$ is divisible.
\end{exa}

\begin{pro}\label{p204}
The following statements are equivalent for an $R$-module $W$:

$(1)$\ \ $W$ is $L_{n}$-injective.

$(2)$\ \ For any $R$-module $M\in \mathcal{F}_n$ and any integer $k \geq 1$, then $\Ext_{R}^{k}(M,W)=0$.

$(3)$\ \ Any exact sequence $0\rightarrow W\rightarrow B \rightarrow C \rightarrow 0$ with $C\in \mathcal{F}_n$ is split.

$(4)$\ \ For any exact sequence $0\rightarrow A \rightarrow B \rightarrow C \rightarrow 0$ with $C\in \mathcal{F}_n$, then
the sequence $0\rightarrow \Hom(C,W) \rightarrow \Hom(B,W) \rightarrow \Hom(A,W) \rightarrow 0$ is exact.
\end{pro}

\no{\bf Proof.}  $(1)\Leftrightarrow (3) \Leftrightarrow (4)$ and $(2)\Rightarrow (1)$ are clear.

$(1)\Rightarrow (2)$. From the definition we have $\Ext_{R}^{1}(M,W)=0$, that is, the assertion is true for the case $k=1$.

Assume $k>1$. Let $0\ra A\ra F\ra M\ra 0$ be exact, where $F$ is projective. Thus $\Ext_{R}^{k}(M,W)\cong \Ext_{R}^{k-1}(A,W)$. Note that $A\in {\mathcal F}_n$.
Hence $\Ext_{R}^{k}(M,W)=0$ by induction on $k$.\hfill$\Box$

\begin{pro}\label{p205} Let $0\rightarrow A \rightarrow B \rightarrow C  \rightarrow 0$ be an exact sequence. Then:

$(1)$\ \ If $A\in\mathcal{L}_n$, then $B\in\mathcal{L}_n$ if and only if $C\in\mathcal{L}_n$.

$(2)$\ \ If $C \in \mathcal{F}_n$, then $A \in \mathcal{F}_n$ if and only if $B \in \mathcal{F}_n$.
\end{pro}

\no{\bf Proof.} These are straightforward.\hfill$\Box$

Let $\mathcal{C}$ be a class of $R$-modules and $M$ a $R$-module. For a $\mathcal{C}$-resolution of $M$
$$\cdots \ra C_n\ra \cdots \rightarrow C_{1} \rightarrow C_{0} \rightarrow M \rightarrow 0, \qquad (\mbox{resp.}\quad
0\rightarrow M \rightarrow C^{0} \rightarrow C^{1} \rightarrow \ra C^n\ra\cdots,)$$
set $K_{0}=M$, $K_{1}=\ker(C_{0} \rightarrow M)$, $K_{i}=\ker(C_{i-1}\rightarrow C_{i-2})$
(resp. $Q^{0}=M$, $Q^{1}=\cok(M \rightarrow C^{0})$, $Q^{i}=\cok(C_{i-2} \rightarrow C_{i-1})$) for $i\geq 2$.
The $n$-th kernel $K_{n}$ (resp. cokernel $Q^{n}$) $(n\geq 0)$ is called the $n$-th $\mathcal{C}$-syzygy (resp. $\mathcal{C}$-cosyzygy) of $M$.
In particular, if $\mathcal{C}$ is the class of projective modules (resp. flat modules), then $K_{n}$ is simply called the $n$-th syzygy
(resp. $n$-york) of $M$; and if $\mathcal{C}$ is the class of injective modules, then $Q^{n}$ is simply called the $n$-th cosyzygy
of $M$.

\begin{thm}\label{t206}
Let $n$ and $m$ be two given nonnegative integers and let $W$ be an $L_{n}$-injective module. Then the
$m$-th cosyzygy of $W$ is $L_{n+m}$-injective.
\end{thm}

\no{\bf Proof.} The case $m=0$ is clear. Now assume $m>0$. Let $W'$ be an $m$-th cosyzygy of $W$ and let $M\in {\mathcal F}_{n+m}$. Let $B$ be an
$(m-1)$-th syzygy of $M$. Then $\fd_{R}B\leq n$. Hence $\Ext_R^1(M,W')\cong \Ext_R^{m+1}(M,W)\cong \Ext_R^1(B,W)=0$. Therefore,
$W'\in {\mathcal L}_{n+m}$.\hfill$\Box$

\begin{cor}\label{c207}
Let $C$ be a cotorsion module. then the $n$-th cosyzygy of $C$ is $L_{n}$-injective.
\end{cor}

Corresponding with the consideration of flat dimension at most $n$, Fuchs and Lee have proved that if $R$ is an integral
domain, then an $R$-module $M$ satisfies $\fd_{R}M \leq n$ if and only if $\Ext_{R}^n(M,W)=0$ for any
$L_1$-injective module $W$ (See \cite[Theorem 2.2]{Le11} and \cite[Lemma 5.5]{FL10}). The following is the further discussion on
the modules with flat dimension $\leq n$.

\begin{thm}\label{t208}
Let $R$ be any ring and $M$ be an $R$-module and $n\geq 1$. Then the following are equivalent:

$(1)$\ \ $\fd_{R}M\leq n$.

$(2)$\ \ $\Ext_R^1(M,N)=0$ for any $L_{n}$-injective module $N$.

$(3)$\ \ $\Ext_{R}^i(M,N)=0$ for any $L_{n}$-injective module $N$ and any $i\geq 1$.

$(4)$\ \ If $0 \leq m \leq n$, then $\Ext_{R}^{n-m+1}(M,N)=0$ for any $L_{m}$-injective module $N$.

$(5)$\ \ If $0 \leq m \leq n$, then $\Ext_{R}^{n-m+i}(M,N)=0$ for any $L_{m}$-injective module $N$ and any $i\geq 1$.

$(6)$\ \ $\Ext_{R}^{n}(M,W)=0$ for any weak-injective {\rm ($L_1$-{\it injective})} module $W$.

$(7)$\ \ $\Ext_{R}^{i}(M,W)=0$ for any weak-injective module $W$ and any $i\geq n$.

$(8)$\ \ $\Ext_{R}^{n+1}(M,C)=0$ for any cotorsion module $C$.

$(9)$\ \ $\Ext_{R}^{n+i}(M,C)=0$ for any cotorsion module $C$ and any $i\geq 1$.
\end{thm}

\no{\bf Proof.}  $(3)\Rightarrow (2)$ and $(9)\Rightarrow (8)$ and $(7) \Rightarrow (6)$ are clear.

(1)$\Ra$(3). It follows by Proposition \ref{p204}.

(5)$\Ra$(9). It follows by picking $m=0$.

(5)$\Ra$(7). It follows by picking $m=1$.

(4)$\Ra$(6). It follows by picking $m=1$.

(5)$\Ra$(3). It follows by picking $m=n$.

(2)$\Ra$(4). Let $0\ra N\ra E_0\ra E_1\ra \cdots\ra E_{n-m-1}\ra W\ra 0$ be exact, where $E_0,E_1,\cdots,E_{n-m-1}$ be injective.
Thus $W$ is an $(n-m)$-cosyzygy of $N$. By Theorem \ref{t206}, $W$ is $L_{n}$-injective. Thus $\Ext_R^{n-m+1}(M,N)\cong\Ext_R^{1}(M,W)=0$.

(3)$\Ra$(5). It is simillar to the proof of $(2)\Rightarrow (4)$.

(6)$\Ra$(8). Let $0\ra C\ra E\ra W\ra 0$ be exact, where $E$ is injective. By Corollary \ref{c207}, $W$ is weak-injective. Hence
$\Ext_{R}^{n+1}(M,C)\cong\Ext_R^n(M,W)=0$.

(8)$\Ra$(1). Let $K$ be an $n$-syzygy of $M$. Then, for any cotorsion module $C$, $\Ext_R^1(K,C)\cong \Ext_R^{n+1}(M,C)=0$.
Hence $K$ is flat by \cite[Lemma 3.4.1]{Xu96}. Therefore, $\fd_RM\leq n$.\hfill$\Box$

Let $R$ be a ring. Bass defined in \cite{Ba60} the left weak finitistic dimension of $R$ as follow:
$$l.\FFD(R)=\sup\{\,\fd_RM\,|\,\mbox{$M$ is an $R$-module with $\fd_RM<\infty$}\,\}.$$

\begin{thm}\label{t211}
Let $n<m$ be two given integers. Then the following statements are equivalent for a ring $R$:

$(1)$\ \ $l.\FFD(R)\leq n$.

$(2)$\ \ ${\mathcal F}_m={\mathcal F}_n$.

$(3)$\ \ $\mathcal{L}_m=\mathcal{L}_n$.

$(4)$\ \ Every $L_n$-injective module is $L_m$-injective.
\end{thm}

\no{\bf Proof.} $(1)\Ra(2)$. It is clear.

$(2)\Ra(3)$. It follows from the facts ${\mathcal L}_m={\mathcal F}_m^{\bot}$ and ${\mathcal L}_n={\mathcal F}_n^{\bot}$.

$(3)\Ra(4)$. It is trivial.

$(4)\Ra(1)$. It is enough that we assume $m=n+1$. Let $M$ be an $R$-module with $\fd_RM=s<\infty$. If $s>n$, without loss of generality, we can assume $s=n+1$.
Hence we have $\Ext_R^1(M,W)=0$ for any $L_n$-injective module $W$ by hypothesis. And so $\fd_RM\leq n$ by
Theorem \ref{t208}, a contradiction. Therefore, $\fd_RM\leq n$, and hence $l.\FFD(R)\leq n$.\hfill$\Box$

\begin{cor}\label{c212}
Let $n\geq 1$. Then every cotorsion $R$-module is $L_{n}$-injective if and only if $l.\FFD(R)=0$.
\end{cor}

\no{\bf Proof.}\ \ Pick $m=0$ in Theorem \ref{t211}.\hfill$\Box$

\vskip5mm
\section{The Cotorsion Theory ($\mathcal{F}_{n}, \mathcal{L}_{n}$)}

\vskip3mm
The goal of this section is to show $(\mathcal{F}_{n}, \mathcal{L}_{n})$ is a complete cotorsion theory for an arbitrary ring $R$.

\begin{defn}\label{d301}
A right $R$-module $D$ is called $L_{n}$-flat $(L$ means Lee$)$ if $\Tor^{R}_{1}(D,M)=0$ for all left $R$-modules $M\in {\mathcal F}_n$.
\end{defn}

Clearly, flat modules are $L_n$-flat. In \cite{EH12} $L_n$-flat modules are called $n$-torsionfree.
In the following we denote ${\mathcal D}_n$ the class of $L_n$-flat modules.

For a right $R$-module $D$, write $D^+=\Hom_Z(D,{\Bbb Q}/{\Bbb Z})$, which is called the character module of $D$.

\begin{lem}\label{l3003}

$(1)$\ \  $M^{+}$ is pure injective for every $R$-module $M$.

$(2)$\ \  An $R$-module $M$ is flat if and only if $M^{+}$ is injective.

$(3)$\ \  Every pure injective $R$-module is cotorsion.

$(4)$\ \  $M^{+}$ is cotorsion for every $R$-module $M$.
\end{lem}

\no{\bf Proof.} $(1)$ and $(2)$ by \cite[Lemma 2.2]{Le06}.

$(3)$ See \cite[Lemma 5.3.23]{EJ00}.

$(4)$ It follows by  $(1)$ and $(3)$.\hfill$\Box$

\begin{pro}\label{p302}
A right $R$-module $D$ is $L_n$-flat if and only if $D^+$ is $L_n$-injective.
\end{pro}

\no{\bf Proof.} It follows by the following standard isomorphism
$$(\Tor^R_1(D,M))^+\cong \Ext_R^1(M,D^+),$$
where $M\in {\mathcal F}_n$.\hfill$\Box$

\begin{lem}\label{l303}
Let $X$ be a right $R$-module. Then any $n$-york of $X$ is $L_n$-flat.
\end{lem}

\no{\bf Proof.} Let $0\ra D\ra D_{n-1}\ra \cdots \ra D_1\ra D_0\ra X\ra 0$ be exact, where $D_0,D_1,\cdots,D_{n-1}$ are flat.
Thus $(D_0)^+,(D_1)^+,\cdots,(D_{n-1})^+$ are injective and $X^+$ is cotorsion by lemma \ref{l3003} and
$0\ra X^+\ra (D_0)^+\ra (D_1)^+\ra\cdots \ra (D_{n-1})^+\ra D^+\ra 0$ is exact. By Corollary \ref{c207}, $D^+$ is $L_n$-injective.
By Proposition \ref{p302}, $D$ is $L_n$-flat.

\begin{thm}\label{t304}
Let $M$ be an $R$-module. Then $\fd_RM\leq n$ if and only if $\Tor^R_1(D,M)=0$ for any $D\in {\mathcal D}_n$.
\end{thm}

\no{\bf Proof.} Suppose $\fd_RM\leq n$. Then it is clear that $\Tor^R_1(D,M)=0$ for any $D\in {\mathcal D}_n$. For the
converse, let $X$ be a right $R$-module and let $D$ be the $n$-th york of $X$. Then $D$ is $L_n$-flat by Lemma \ref{l303}. Therefore,
$\Tor^R_{n+1}(X,M)\cong \Tor^R_1(Y,M)=0$ by hypothesis. Hence $\fd_RM\leq n$.\hfill$\Box$

Let $(\mathcal{A,B})$ be a cotorsion theory. Recall that $(\mathcal{A,B})$ is called hereditary if whenever
$0\rightarrow A \rightarrow B \rightarrow C \rightarrow 0$ is an exact sequence with $B,C\in \mathcal{A}$, then $A$ is also in $\mathcal{A}$;
equivalently, whenever $0\rightarrow A \rightarrow B \rightarrow C \rightarrow 0$ is an exact sequence with $A,B\in \mathcal{B}$, then $C$
is also in $\mathcal{B}$. Recall that $(\mathcal{A,B})$ is said to be complete \cite{Tr00} if every $R$-module has a special $\mathcal{A}$-precover.
By \cite[Lemma 1.13]{Tr00}, a cotorsion theory $(\mathcal{A,B})$ is complete if and only if every $R$-module has a special $\mathcal{B}$-preenvelope.

Let ${\mathcal A}$ be a class of left $R$-modules and let ${\mathcal B}$ be a class of right $R$-modules. Write
$${}^{\top}{\mathcal A}=\{D\in {\frak M}_R\,|\,\mbox{$\Tor^R_1(D,M)=0$ for any $M\in {\mathcal A}$}\},$$
and
$${\mathcal B}^{\top}=\{Y\in {}_R{\frak M}\,|\,\mbox{$\Tor^R_1(D,Y)=0$ for any $D\in {\mathcal B}$}\}.$$
If ${\mathcal B}={}^{\top}{\mathcal A}$ and ${\mathcal A}={\mathcal B}^{\top}$, then the pair $({\mathcal A},{\mathcal B})$ is called a
Tor-torsion theory.

For a class ${\mathcal C}$ of modules, set
$${\frak S}_{\mathcal C}=({}^{\bot}{\mathcal C},({}^{\bot}{\mathcal C})^{\bot}).$$

\begin{lem}\label{l305}
$(1)$\ \ If $({\mathcal A},{\mathcal B})$ be a Tor-torsion theory, then ${\frak C}:=({\mathcal A},{\mathcal A}^{\bot})$ is a cotorsion theory.
Moreover, if we write ${\mathcal C}=\{B^+\,|\,B\in {\mathcal B}\}$, then ${\mathcal C}$ is a subclass of the class of pure-injective modules
and ${\frak S}_{\mathcal C}={\frak C}$.

$(2)$\ \ Let ${\frak S}_{\mathcal C}=({\mathcal A},{\mathcal B})$ be a cotorsion theory. If ${\mathcal C}$ is a subclass of the class of pure-injective modules,
the ${\frak S}_{\mathcal C}$ is a complete cotorsion theory. Moreover, every $R$-module has an ${\mathcal A}$-cover and a ${\mathcal B}$-envelope.
\end{lem}

\no{\bf Proof.}\ \ See \cite[Lemma 1.11\& Theorem 2.8]{Tr00}.\hfill$\Box$

\begin{thm}\label{t306}
$(1)$\ \ $({\mathcal F}_n,{\mathcal D}_n)$ is a Tor-torsion theory.

$(2)$\ \ $(\mathcal{F}_{n},\mathcal{L}_{n})$ is a complete cotorsion theory.

$(3)$\ \ Every $R$-module has a special ${\mathcal F}_n$-cover and a special ${\mathcal L}_n$-envelope. Further, for any $R$-modulea $M$ and $N$,
there are exact sequences
$$0\ra A\ra F\ra M\ra 0\qquad \mbox{and} \qquad 0\ra N\ra W\ra B\ra 0,\eqno{(3.1)}$$
where $F$ is the ${\mathcal F}_n$-cover, $W$ is the ${\mathcal L}_n$-envelope, $B\in {\mathcal F}_n$, and $A\in {\mathcal L}_n$.
\end{thm}

\no{\bf Proof.} (1)\ \ Clearly, ${\mathcal D}_n={}^{\top}{\mathcal F}_n$. By Theorem \ref{t304}, ${\mathcal F}_n={\mathcal D}_n^{\top}$.

(2)\ \ Clearly, ${\mathcal L}_n={\mathcal F}_n^{\bot}$. Hence $(\mathcal{F}_{n},\mathcal{L}_{n})$ is a cotorsion theory by Lemma \ref{l305} (1).
Moreover, $({\mathcal F}_n,{\mathcal L}_n)$ is a complete cotorsion theory by Lemma \ref{l305} (2).

(3)\ \ By Lemma \ref{l305} (2), $M$ has the ${\mathcal F}_n$-cover $\phi:F\ra M$ and $N$ has the ${\mathcal L}_n$-envelope $\varphi:N\ra W$. Since
$M$ has a a special ${\mathcal F}_n$-precover, $\phi$ is epic. Because ${\mathcal F}_n$ is closed under extension, $A:=\ker(\phi)\in {\mathcal L}_n$
by \cite[Lemma 2.1.1]{Xu96}. Hence $\phi:F\ra M$ is also special.

The other statement is dual to the argument above, but we need to apply \cite[Lemma 2.1.2]{Xu96}.\hfill$\Box$

In the following we denote $F_n(M)$ the ${\mathcal F}_n$-cover of an $R$-module $M$ and $L_n(N)$ the ${\mathcal L}_n$-envelope of
an $R$-module $N$.

\begin{thm}\label{t307}
Let $M$ and $N$ be $R$-modules. Then the following statements are true:

$(1)$\ \ $M$ is $L_{n}$-injective if and only if $F_n(M)$ is $L_{n}$-injective.

$(2)$\ \ $\fd_{R}N\leq n$ if and only if $\fd_{R}L_n(N)\leq n$.
\end{thm}

\no{\bf Proof.} It follows directly from the exact sequences in (3.1)\hfill$\Box$

\begin{cor}\label{c308}
Let $M$ and $N$ be $R$-modules. Then the following statements are true:

$(1)$\ \ $M$ is cotorsion if and only if so is its flat cover.

$(2)$\ \ $N$ is flat if and only if so is its cotorsion envelope.
\end{cor}

\begin{thm}\label{tshg45s}

Let $N\suse C$ be an extension of $R$-modules and $C$ be $L_n$-injective. Then the following statements are equivalent:

$(1)$\ \ $C$ is an ${\mathcal L}_n$-envelope of $N$.

$(2)$\ \ $C/N\in {\mathcal F}_n$, and there is no submodule $0\neq A \suse C$ such that $N\bigcap A=0$ and $C/(N+A)\in {\mathcal F}_n$.
\end{thm}

\no{\bf Proof.}
(1)$\Ra$(2). Let $C$ be an ${\mathcal L}_n$-envelope of $N$ and $i:N\ra C$ be an inclusion homomorphism.
$C/N\in {\mathcal F}_n$ holds by Theorem \ref{t306}. Now, let $A$ be a submodule of $C$ such that $A\bigcap N=0$
and $B:=C/(N+A)\in {\mathcal F}_n$. Then the natural homomorphism $\phi:N\ra C/A$ is monic and $\cok(\phi)=C/(N+A)$. Let
$L$ be an ${\mathcal L}_n$-envelope of $C/A$. Then we get the following commutative diagram with exact rows
$$\xymatrix@R=19pt@C=27pt{
&&0\ar[d]&0\ar[d]\\
0\ar[r]&N\ar[r]^{\phi}\ar@{=}[d]&C/A\ar[r]\ar[d]&B\ar[r]\ar[d]&0\\
0\ar[r]&N\ar[r]^{\phi}&L\ar[r]\ar[d]&X\ar[r]\ar[d]&0\\
&&Y\ar@{=}[r]\ar[d]&Y\ar[d]\\
&&0&0}$$
$X\in {\mathcal F}_n$ since $Y,B\in {\mathcal F}_n$. Let $\pi:C\ra C/A$ be a natural homomorphism
and $\lambda:C/A\ra L$ an inclusion homomorphism. Let $f:L\ra C$ be a homomorphism  such that
$f\phi=i$. Then $f\lambda\pi$ is an isomorphism and $\pi$ is monic. Hence $A=0$.

(2)$\Ra$(1).  Let $E$ be an ${\mathcal L}_n$-envelope of $N$ and $\lambda:N\ra E$ an inclusion homomorphism.
Then $D:=E/N\in {\mathcal F}_n$. So there exist homomorphisms $f:E\ra C, g:C\ra E$  such that $f\lambda=i, gi=\lambda$. Thus
$(gf)\lambda=\lambda$, where $gf:E\ra E$. So
$fg$ is an isomorphism and there exist homomorphism $h:E\ra E$ such that $hgf={\bf 1}_E$. Thus $C=\im(f)\bigoplus \ker(hg)$. Set
$A=\ker(hg)$. Then $N\bigcap A=0$ and $C/(N+A)\cong E/N\in {\mathcal F}_n$.
By hypothesis, $A=0$.  Hen $f:E\ra C$ is an isomorphism.
\hfill$\Box$

\begin{thm}\label{t309}
Let $n<m$ be two given integers. Then the following statements are equivalent for a ring $R$:

$(1)$\ \ $l.\FFD(R)\leq n$.

$(2)$\ \ ${\mathcal D}_m={\mathcal D}_n$.

$(3)$\ \ Every $L_n$-flat module is $L_m$-flat.
\end{thm}

\no{\bf Proof.} (1)$\Ra$(2). By Theorem \ref{t211}, ${\mathcal F}_m={\mathcal F}_n$. Hence, by Theorem \ref{t306},
$${\mathcal D}_m={}^{\top}{\mathcal F}_m={}^{\top}{\mathcal F}_n={\mathcal D}_n.$$

(2)$\Ra$(1). ${\mathcal F}_m={\mathcal D}_m^{\top}={\mathcal D}_n^{\top}={\mathcal F}_n$ by applying Theorem \ref{t306} again.

(2)$\lra$(3). Clearly.\hfill$\Box$

\vskip5mm
\section{$L_{n}$-Injective dimensions of modules}

In this section, we study the $L_{n}$-injective modules over an arbitrary ring $R$.

Let $M$ be an $R$-module. If there exsits an exact sequence
$$0\rightarrow M \rightarrow W_{0} \rightarrow W_{1} \rightarrow \cdots \rightarrow W_{m-1} \rightarrow W_{m} \rightarrow \cdots \qquad \eqno{(4.1)}$$
in which each $W_{i}$ is $L_{n}$-injective, then this exact sequence is called an $L_{n}$-injective resolution of $M$. Certainly, every $R$-module $M$
has an $L_n$-injective resolution. If the following homomorphisms
$$M \rightarrow W_{0},\quad \cok(M \rightarrow W_{0})\rightarrow W_{1}, \quad  \cdots, \quad \cok(W_{i-2} \rightarrow W_{i-1})\rightarrow W_{i}, \cdots.
\eqno{(4.2)}$$
are ${\mathcal L}_n$-envelopes, then the sequence (4.1) is called a minimal $L_{n}$-injective resolution of $M$.

\begin{pro}\label{p401}
Every $R$-module $M$ has a minimal $L_n$-injective resolution.
\end{pro}

\no{\bf Proof.} By Theorem \ref{t306}, $M$ has an $\mathcal{L}_{n}$-envelope $W_{0}$. Note $C_{0}=\cok(M \rightarrow W_{0})$.
Then $C_{0}$ also has an $\mathcal{L}_{n}$-envelope $W_{1}$ with $C_{1}=\cok(C_{0} \rightarrow W_{1})$ also by Theorem \ref{t306}.
The resuit holds by repeating this process.\hfill$\Box$

\begin{defn}\label{d402}
Let $M$ be an $R$-module. By the $\mathcal{L}_{n}$-injective dimension $L_nid_{R}M$ of $M$ is defined to be the smallest integer $m\geq 0$ such that the sequence
$0\rightarrow M \rightarrow W_{0} \rightarrow W_{1} \rightarrow \cdots \rightarrow W_{m-1} \rightarrow W_{m} \rightarrow 0$ in which
each $W_{i}$ is $L_{n}$-injective for $0\leq i \leq m$ is exact. If there is no such integer $m$, set $L_nid_{R}M=\infty$.
\end{defn}

\begin{exa}\label{e403}
Let $M$ be an $R$-module.

$(1)$\ \ $M$ is $L_{n}$-injective if and only if $\lid_{R}M=0$.

$(2)$\ \ $\lid_{R}M\leq id_{R}M$.

$(3)$\ \ If $m\leq n$, then $L_m\id_{R}M\leq \lid_{R}M$ since every $L_{n}$-injective module is $L_{m}$-injective.
\end{exa}

\begin{thm}\label{t404}
Let $m$ be a nonnegative integer. The following statements are equivalent for an $R$-module $N$:

$(1)$\ \ $\lid_{R}N\leq m$.

$(2)$\ \ $\Ext_{R}^{m+1}(M,N)=0$ for any $M\in \mathcal{F}_{n}$.

$(3)$\ \ $\Ext_{R}^{m+i}(M,N)=0$ for any $M\in \mathcal{F}_{n}$ and any $i\geq 1$.

$(4)$\ \ If $0\rightarrow N\rightarrow W_{0} \rightarrow W_{1} \rightarrow \cdots \rightarrow W_{m-1} \rightarrow W_{m} \ra 0$
is exact, where $W_0,W_1,\cdots,W_{m-1}$ are $L_{n}$-injective, then $W_m$ is $L_{n}$-injective.

$(5)$\ \ If $0\rightarrow N\rightarrow W_{0} \rightarrow W_{1} \rightarrow \cdots \rightarrow W_{m-1} \rightarrow W_{m} \rightarrow 0$ is
exact, where $W_0,W_1,\cdots,W_{m-1}$ is injective, then $W_{m}$ is $L_{n}$-injective.

$(6)$\ \ The $m$-cosyzygy of $N$ in its minimal $L_n$-injective resolution is $L_{n}$-injective.
\end{thm}

\no{\bf Proof.}  (3)$\Ra$(4)$\Ra$(5)$\Ra$(6). Trivially.

(1)$\Ra$(2). Since $L_n\id_RN\leq m$, there is an exact sequence
$0\rightarrow N \rightarrow W_{0} \rightarrow W_{1} \rightarrow \cdots \rightarrow W_{m-1} \rightarrow W_{m}\rightarrow 0$ in which
each $W_{i}$ is $L_{n}$-injective. Therefore, for $M\in {\mathcal F}_n$, $\Ext_R^{m+1}(M,N)\cong \Ext_R^1(M,W_m)=0$.

(2)$\Ra$(3). For $M\in \mathcal{F}_{n}$, there is an exact sequence $0 \rightarrow K \rightarrow P \rightarrow M\rightarrow 0$
with $P$ projective. Then $\Ext_{R}^k(M,N)\cong \Ext_R^{k-1}(K,N)$. Thus it follows by induction.

(6)$\Ra$(1). Take any minimal $L_n$-injective resolution of $M: 0\rightarrow M \rightarrow E_{0} \rightarrow E_{1} \rightarrow \cdots \rightarrow E_{m-1} \rightarrow L^{m-1} \rightarrow 0$ with
each $E_{i}$ be $L_n$-injective and $L^{m-1}$ be the $m$-th cosyzygy of $M$. Then the result follows since $L^{m-1}$ is $L_{n}$-injective by hypothesis.  \hfill$\Box$

\begin{pro}\label{p405}
Let $0\rightarrow A \rightarrow W\rightarrow C \rightarrow 0$ be an exact sequence, where $W$ is $L_{n}$-injective.

$(1)$\ \ If $A$ is $L_n$-injective, then so is $C$.

$(2)$\ \ If $\lid_RA>0$, then $\lid_{R}C=\lid_{R}A-1$.
\end{pro}

\no{\bf Proof.} (1)\ \ It follows by Theorem \ref{t206}.

(2)\ \ It follows from the isomorphism $\Ext_R^m(M,C)\cong \Ext_R^{m+1}(M,A)$, where $M\in{\mathcal F}_n$.\hfill$\Box$

\begin{pro}\label{p406}
Let $0\rightarrow A \rightarrow B\rightarrow C \rightarrow 0$ be an exact sequence.

$(1)$\ \ If two of $\lid_{R}A$, $\lid_{R}B$, and $\lid_{R}C$ are finite, so is the third.

$(2)$\ \ If $A$ is $L_n$-injective, then $\lid_RB=\lid_RC$.

$(3)$\ \ $\lid_{R}B\leq \max\{\lid_{R}A,\lid_{R}C\}$.

$(4)$\ \ $\lid_{R}C\leq \max\{\lid_{R}A-1,\lid_{R}B\}$.

$(5)$\ \ $\lid_{R}A\leq \max\{\lid_{R}B,\lid_{R}C+1\}$.
\end{pro}

\no{\bf Proof.} These results follow easily from the exact sequence
$$\Ext_R^m(M,C)\ra \Ext_R^{m+1}(M,A)\ra \Ext_R^{m+1}(M,B)\ra \Ext_R^{m+1}(M,C)\ra\Ext_R^{m+2}(M,A)$$
and applying Theorem \ref{t404}, where $M\in{\mathcal F}_n$.\hfill$\Box$

\begin{pro}\label{p407}
Let $\{M_i\,|\,i\in\Gamma\,\}$ be a family of $R$-modules. Then
$$\lid_{R}(\mbox{$\prod\limits_{i\in\Gamma}M_i$})=\sup\{\lid_{R}M_i\,|\,i\in \Gamma\,\}.$$
\end{pro}

\no{\bf Proof.} It is straightforward.\hfill$\Box$

\begin{thm}\label{t408}
Let $N$ be an $R$-module. Then $\lid_{R}F_n(N)=\lid_{R}N$.
\end{thm}

\no{\bf Proof.} This follows from the first exact sequence in Theorem \ref{t306} (3) and Proposition \ref{p406}.\hfill$\Box$

\begin{pro}\label{p409}
Let $N$ be an $R$-module. If $id_RN<\infty$ and every injective $R$-module has the flat dimension at most $n$,
then $\lid_{R}N=\id_RN$.
\end{pro}

\no{\bf Proof.} Write $\id_RN=m$. Then $\lid_{R}N\leq m$. Pick an injective $R$-module $E$ such that $\Ext_{R}^{m}(E,N)\neq 0$.
As $E\in {\mathcal F}_n$ we get $\lid_{R}N\geq m$. Hence $\lid_RN=m$.\hfill$\Box$

\begin{thm}\label{t410}
Let $N$ be an $R$-module with $\lid_{R}N=m<\infty$. Then there exsits an $L_n$-injective module $W\in {\mathcal F}_n$
such that $\Ext_{R}^{m}(W,N)\neq 0$.
\end{thm}

\no{\bf Proof.} Since $\lid_{R}N=m$, we pick an $R$-module $M\in \mathcal{F}_{n}$ such that $\Ext_{R}^{m}(M,N)\neq 0$.
Let $0\ra M\ra W\ra B\ra 0$ be exact, where $W$ is the $L_n$-injective envelope. By Theorem \ref{t306} $B\in {\mathcal F}_n$.
By Theorem \ref{t307}, $W\in {\mathcal F}_n$. Since the sequence $\Ext_R^m(W,N)\ra \Ext_R^m(M,N)\ra \Ext_R^{m+1}(B,N)=0$
is exact, we obtain $\Ext_{R}^{m}(W,N)\neq 0$.\hfill$\Box$

\vskip5mm
\section{$L_n$-Global dimensions of a ring}

\vskip3mm
To characterize propoties of rings by using $L_n$-injectivity, we are in the position to define the $L_n$-global dimension of
a ring.

\begin{defn}\label{d50l}
For a ring $R$, its left $L_n$-global dimension $l.\ld(R)$ is defined by
$$l.\ld(R)=\sup\{\lid_{R}M\,|\,\mbox{$M$ is an $R$-module.}\}.$$
\end{defn}

\begin{exa}\label{e502}
For a ring $R$, we have:

$(1)$\ \ $l.\ld(R)\leq l.\gl(R)$.

$(2)$\ \ If $m\leq n$, then $l.L_m\dim(R)\leq l.\ld(R)$.
\end{exa}

\begin{exa}\label{e503}
By {\rm \cite[Lemma 3.6]{Le06}} and {\rm \cite[Corollary 6.4]{FL09}}, a domain $R$ is APD if and only if $L_1\dim(R)\leq 1$.
\end{exa}

\begin{thm}\label{t504}
Let $m$ be a nonnegative integer. Then the following statements are equivalent for a ring $R$:

$(1)$\ \ $l.\ld(R)\leq m$.

$(2)$\ \ $\Ext_{R}^{m+i}(M,N)=0$ for any $M\in \mathcal{F}_{n}$ and $N\in {}_R{\frak M}$ and for any $i\geq 1$.

$(3)$\ \ $\Ext_{R}^{m+1}(M,N)=0$ for any $M\in \mathcal{F}_{n}$ and $N\in {}_R{\frak M}$.

$(4)$\ \ $\Ext_{R}^{m+i}(M,N)=0$ for any $M,N\in \mathcal{F}_{n}$ and any $i\geq 1$.

$(5)$\ \ $\Ext_{R}^{m+1}(M,N)=0$ for any $M,N\in \mathcal{F}_{n}$.

$(6)$\ \ $\sup\{\lid_{R}N\,|\,N\in \mathcal{F}_{n}\}\leq m$.

$(7)$\ \ $\sup\{\pd_{R}M\,|\,M\in \mathcal{F}_{n}\}\leq m$.
\end{thm}

\no{\bf Proof.} (2)$\Ra$(4)$\Ra$(5) and (3)$\Ra$(7) are trivial.

(1)$\lra$(2)$\lra$(3). It follows from Theorem \ref{t404}.

(5)$\Ra$(6). Let $N\in {\mathcal F}_n$. Then $\lid_RN\leq m$ by applying Theorem \ref{t404} again. Therefore, $\sup\{\lid_{R}N\,|\,N\in \mathcal{F}_{n}\}\leq m$.

(6)$\Ra$(1). Let $N\in {}_R{\frak M}$. Then we have the exact sequence $0\ra A\ra F_n(N)\ra N\ra 0$ by Theorem \ref{t306},
where $A\in {\mathcal L}_n$. By Proposition \ref{p406}, $\lid_RN=\lid_RF_n(N)\leq m$. Therefore,
$l.\ld(R)\leq m$.

(7)$\Ra$(3). Let $M\in {\mathcal F}_n$ and $N\in {}_R{\frak M}$. Since $\pd_RM\leq m$, we have $\Ext_R^{m+1}(M,N)=0$. \hfill$\Box$

\begin{cor}\label{c505}
For any ring $R$, the following are identical:

$(1)$\ \ $l.\ld(R)$.

$(2)$\ \ $\sup\{\lid_{R}N\,|\,N\in \mathcal{F}_{n}\}$.

$(3)$\ \ $\sup\{\pd_{R}M\,|\,M\in \mathcal{F}_{n}\}$.
\end{cor}

\begin{thm}\label{t506}
Let $m$ be a nonnegative integer. If $l.\ld(R)<\infty$, then the following statements are equivalent:

$(1)$\ \ $l.\ld(R)\leq m$.

$(2)$\ \ $\sup\{\pd_{R}M\,|\,M\in {\mathcal F}_{n}\bigcap {\mathcal L}_n\}\leq m$.

$(3)$\ \ $\sup\{\pd_{R}L_n(M)\,|\,M\in \mathcal{F}_{n}\}\leq m$.

$(4)$\ \ $\sup\{\pd_{R}F_n(M)\,|\,M\in {\mathcal L}_{n}\}\leq m$.

$(5)$\ \ $\sup\{\pd_{R}F_n(M)\,|\,M\in {}_R{\frak M}\}\leq m$.

$(6)$\ \ $\sup\{\lid_{R}M\,|\,\mbox{$M$ is projective}\}\leq m$.
\end{thm}

\no{\bf Proof.} (1)$\Ra$(5). It follows from Theorem \ref{t504}.

(5)$\Ra$(4). It is trivial.

(4)$\Ra$(2). For any $M\in {\mathcal F}_{n}\bigcap {\mathcal L}_n$, $F_n(M)=M$, and hence $\pd_RM\lst m$ by hypothesis.

(2)$\Ra$(1). Let $M\in {\mathcal F}_n$. Because $l.\ld(R)<\infty$, by Proposition \ref{p407} and Theorem \ref{t306}
there is an exact sequence
$$0\ra M\ra W_0\ra W_1\ra\cdots \ra W_{k-1}\ra W_k\ra 0\eqno{(5.1)}$$
such that $W_0,W_1,\cdots,W_{k-1}$ are $L_n$-envelopes of some modules with the flat dimension at most $n$, and
every ${\mathcal L}_{n}$-cosyzygy
in (5.1) is in ${\mathcal F}_n$. Thus $W_k\in {\mathcal F}_n\bigcap {\mathcal L}_n$. By hypothesis and Theorem \ref{t307},
$\pd_RW_i\leq m$ for $i=0,1,\cdots,k$. These imply that $\pd_RM\leq m$. Consequently, $l.\ld(R)\leq m$ by Theorem \ref{t504}.

(2)$\Ra$(3). By Theorem \ref{t307}, $L_n(M)\in {\mathcal F}_n\bigcap {\mathcal L}_n$.

(3)$\Ra$(2). For any $M\in {\mathcal F}_{n}\bigcap {\mathcal L}_n$, $L_n(M)=M$, and hence $\pd_RM\lst m$ by hypothesis.

(1)$\Ra$(6). It is clear because every projective module is in ${\mathcal F}_n$.

(6)$\Ra$(1). Let $N\in \mathcal{F}_{n}$. Since $l.\ld(R)<\infty$, $\pd_{R}N<\infty$ by applying Theorem \ref{t504}. Then we can pick
a projective resolution of
$$0\rightarrow P_{k} \rightarrow P_{k-1}\rightarrow\cdots \rightarrow P_{1}\rightarrow P_{0}\rightarrow N\rightarrow 0.\eqno{(5.2)}$$
Because $\lid_RP_i\leq m$ by hypothesis for $i=0,1,\cdots,k$, $\lid_RN\leq m$ by applying repeatedly Proposition \ref{p406}.
So $l.\ld(R)\leq m$ by Theorem \ref{t504}.\hfill$\Box$

\begin{cor}\label{c507}
Let $n\leq m$. The following statements are equivalent for any ring $R$:

$(1)$\ \ $l.\ld(R)\leq m$.

$(2)$\ \ If $\fd_RM\leq n$, then $\pd_RM\leq m$.

$(3)$\ \ If $M$ is a submodule of a projective module $P$ with $M\in \mathcal{F}_{n-1}$, then $\pd_{R}M\leq m-1$.
\end{cor}

\no{\bf Proof.}\ \ It is easy by Theorem \ref{t504}.\hfill$\Box$

\begin{cor}\label{c508}
$l.\ld(R)\leq n$ if and only if $\pd_RM\leq n$ for any $M\in{\mathcal F}_n$.
\end{cor}

\no{\bf Proof.}\ \ It follows by taking $m=n$ in Corollary \ref{c507}.\hfill$\Box$

Recall that a ring $R$ is called left $m$-perfect if the projective dimension of every flat module is at most $m$ (see \cite{EJL05}).
By the left global cotorsion dimension ($l.cot.D(R)$) introduced by Mao and Ding \cite{MD062} we have:

\begin{cor}\label{c509}
{\rm\cite[Corollary 19.27]{MD062}} A ring $R$ is left $m$-perfect if and only if $l.cot.D(R)\leq m$.
\end{cor}

\no{\bf Proof.}\ \ It follows by taking $n=0$ in Corollary \ref{c507}.\hfill$\Box$

\begin{thm}\label{t510}
Let $m$ and $k$ be nonnegative integers with $n\leq k$ and let $l.\ld(R)\leq m$. If $M\in \mathcal{F}_{k}$, then $\pd_{R}M\leq m+k-n$.
\end{thm}

\no{\bf Proof.}\ \ Let $K$ be the $(k-n)$-th syzygy of $M$ and let $N$ be any $R$-module. Then $\fd_{R}(K)\leq n$.
Hence $\Ext_{R}^{m+k-n+1}(M,N)\cong\Ext_{R}^{m+1}(K,N)=0$ by Theorem \ref{t504}. So $\pd_{R}(M)\leq m+k-n$.\hfill$\Box$

\begin{cor}\label{c511}
Let $k\geq n$. If $\wgl(R)\leq k$ and $l.\ld(R)\leq m$, then $\gl(R)\leq m+k-n$.
\end{cor}

\vskip5mm
\section{The characterizations of rings}

\vskip3mm
In this section we decide first when every $L_n$-injective module is injective.

\begin{thm}\label{t601}
The following statements are equivalent for a ring $R$:

$(1)$\ \ $\wgl(R) \leq n$.

$(2)$\ \ Every $L_{n}$-injective module is injective.

$(3)$\ \ If $N\in {\mathcal L}_n$, then $\fd_{R}N\leq n$.

$(4)$\ \ $\Ext_{R}^{1}(M,W)=0$ for any $M,W\in\mathcal{L}_{n}$.

$(5)$\ \ $\Ext_{R}^{i}(M,W)=0$ for any $M,W\in\mathcal{L}_{n}$ and any $i\geq 1$.

$(6)$\ \ $\fd_{R}L_n(M)\leq n$ for any $M\in{}_R{\frak M}$.

\end{thm}

\no{\bf Proof.} (2)$\Ra$(5)$\Ra$(4) and (1)$\Ra$(6) are trivial.

(4)$\Ra$(3). By Theorem \ref{t208}.

(6)$\Ra$(1). By Theorem \ref{t307}, $\fd_RM\leq n$. Hence $\wgl(R)\leq n$.

(1)$\Ra$(2). Let $M$ be any $R$-module and let $W$ be an $L_n$-injective module. Then $\fd_RM\leq n$ by
hypothesis. Hence we have $\Ext_{R}^{1}(M,W)=0$. Consequently, $W$ is injective.

(2)$\Ra$(1). Let $M$ be any $R$-module and let $W$ be any $L_n$-injective module. Because $W$ is injective by hypothesis,
$\Ext_R^1(M,W)=0$. Therefore, $\fd_RM\leq n$ follows by Theorem \ref{t306}, and hence $\wgl(R)\leq n$.

(3)$\Ra$(2). Let $W$ be an $L_n$-injective module and let $0\ra W\ra E\ra C\ra 0$ is an exact sequence, where $E$ is injective.
Then $C$ is also $L_n$-injective. By hypothesis, $\fd_RC\leq n$. Thus the exact sequence is split. Consequently, $W$ is injective.
\hfill$\Box$

\begin{thm}\label{t602}
Let $\wgl(R)< \infty$. Then the following statements are equivalent:

$(1)$\ \ $\wgl(R) \leq n$.

$(2)$\ \ $F_n(W)$ is injective for any $W\in {\mathcal L}_n$.

$(3)$\ \ If $W\in {\mathcal L}_n\bigcap {\mathcal F}_n$, then $W$ is injective.

$(4)$\ \ $L_n(N)$ is injective for any $N\in\mathcal{F}_{n}$.
\end{thm}

\no{\bf Proof.} (1)$\Ra$(2). By Theorem \ref{t307}, $F_n(W)$ is $L_n$-injective. By Theorem \ref{t601}, $F_n(W)$ is
injective.

(2)$\Ra$(3). It is trival as $F_n(W)=W$.

(3)$\Ra$(1). If there is an $R$-module with the flat dimension large than $n$, then there is an $R$-module
$M$ with $\fd_RM=n+1$. Take the exact $0\ra A\ra F_n(M)\ra M\ra 0$ as in (3.1). Then $A$ is $L_n$-injective with $\fd_RA=n$.
Hence $A$ is injective by hypothesis. Therefore, the given sequence is split. Hence $\fd_RM=n$, a contradiction.
Consequently, $\wgl(R)\leq n$.

(3)$\Ra$(4). By Theorem \ref{t307}, $L_n(N)\in {\mathcal L}_n\bigcap {\mathcal F}_n$. Hence $L_n(N)$ is injective by hypothesis.
If $k>n$,

(4)$\Ra$(3). It is trival as $L_n(W)=W$.\hfill$\Box$

\begin{cor}\label{c603}
The following statements are equivalent for a ring $R$:

$(1)$\ \ $R$ is a Von Neumann Regular ring.

$(2)$\ \ Every cotorsion module $N$ is injective.

$(3)$\ \ Every cotorsion module $N$ is flat.

$(4)$\ \ $\Ext_{R}^{1}(M,N)=0$ for any cotorsion modules $M,N$.

$(5)$\ \ $\Ext_{R}^{i}(M,N)=0$ for any any cotorsion modules $M,N$ and any $i\geq 1$.

$(6)$\ \ The flat covers of any cotorsion module $M$ are injective.
\end{cor}

\no{\bf Proof.} Take $n=0$ in Theorem \ref{t601}.\hfill$\Box$

Let $R$ be a ring. Bass defined in \cite{Ba60} the left finitistic projective dimension and left finitistic
flat dimension of $R$ as follow:
$$l.\FPD(R)=\sup\{\pd_RM\,|\,\pd_RM<\infty\}.$$
$$l.\FFD(R)=\sup\{\fd_RM\,|\,\fd_RM<\infty\}.$$

\begin{thm}\label{t604}
For any ring $R$, $l.\ld(R)\leq l.\FPD(R)$. Moreover, if $m=l.\FPD(R)<\infty$, then $l.\FPD(R)=l.\ld(R)$ for ang $n>m$.
\end{thm}

\no{\bf Proof.}\ \ Let $l.\FPD(R)=k<\infty$ and let $M\in \mathcal{F}_{n}$. By \cite[Proposition 6]{Je70}, $\pd_{R}M\leq k$,
and hence $\Ext_{R}^{k+1}(M,N)=0$. Therefore, $l.\ld(R)\leq k$.

If $m=l.\FPD(R)<\infty$, there exist an $R$-module $M$ with $m=\pd_{R}M<n$. Then there exist an $R$-module $N$
such that $\Ext_{R}^{i}(M,N)\neq 0$. Hence $\lid_{R}N\geq m$ since $M\in  \mathcal{F}_{n}$, and $l.\ld(R)\geq m$. That is, $l.\ld(R)=l.\FPD(R)$.

\hfill$\Box$

A ring $R$ is called left perfect if every $R$-module has projective cover, equivaently. every flat $R$-module is projective.
If $R$ is commutative ring, then $R$ is perfect if and only $\FPD(R)=0$ (see \cite{Va76}). By using the
notion of $L_n$-injective modules, we can characterize left perfect rings.

\begin{thm}\label{t605}
Let $n\geq 1$. The following statementes are equivalent for a ring $R$:

$(1)$\ \ $l.\ld(R)=0$, that is, every module is $L_{n}$-injective.

$(2)$\ \ $M$ is projective for any $M \in \mathcal{F}_{n}$.

$(3)$\ \ $M$ is $L_{n}$-injective for any $M\in\mathcal{F}_{n}$.

$(4)$\ \ $\Ext_{R}^{1}(M,N)=0$ for any $M,N\in \mathcal{F}_{n}$.

$(5)$\ \ $\Ext_{R}^{i}(M,N)=0$ for any $M,N\in \mathcal{F}_{n}$ and $i\geq 1$.

$(6)$\ \ $L_n(M)$ is projective for any $M\in \mathcal{F}_{n}$.

$(7)$\ \ $F_{n}(M)$ is projective for any $M\in \mathcal{L}_{n}$.

$(8)$\ \ $F_n(M)$ is $L_{n}$-injective for any $M\in {}_R{\frak M}$.

$(9)$\ \ $R$ is left perfect and $l.\FFD(R)=0$.

$(10)$\ \ $l.\FPD(R)=0$.
\end{thm}

\no{\bf Proof.} $(1)\lra (2)\lra (3)\lra (4)\lra (5)\lra (6)\lra (7)$. By Theorem \ref{t504}.

(1)$\Ra$(8) and $(9)\lra (10)$. Trivial.

(8)$\Ra$(1). By Theorem \ref{t408}.

(1)$\Ra$(9). Since every $L_n$-injective is cotorsion, we have that $R$ is left perfect assertion is true by \cite[Proposition 3.3.1]{Xu96}.
By Corollory \ref{c212}, $l.\FFD(R)=0$.

(9)$\Ra$(1). By Applying \cite[Proposition 3.3.1]{Xu96} and Corollary \ref{c212} again.\hfill$\Box$

\begin{cor}\label{c606}
Let $n\geq 1$.

$(1)$\ \ If $l.\FPD(R)=0$, then every $R$-module is $L_n$-injective.

$(2)$\ \ A commutative $R$ is perfect if and only if every $R$-module is $L_n$-injective.
\end{cor}

\no{\bf Proof.} (1)\ \ It follows by Theorem \ref{t604}.

(2)\ \ It is immediate from (1).\hfill$\Box$

\begin{defn}\label{d607}
A ring $R$ is called left $L_n$-hereditary if every quotient module of an $L_n$-injective
module is $L_n$-injective. In other words, $l.\ld(R)\leq 1$.
\end{defn}

\vskip2mm
\begin{exa}\label{e608}
$(1)$\ \ Left hereditary rings are certainly $L_n$-hereditary for all $n\geq 0$.

$(2)$\ \ Recall that a commtaive ring $R$ is called almost perfect if its proper epic images are perfect. An almost perfect domain is said simply
an APD. In {\rm \cite{Sa10}} it is shown that if $R$ is almost perfect, then $R$ is either perfect or APD.
By {\rm \cite[Lemma 3.6]{Le06}} and {\rm \cite[Corollary 6.4]{FL09}}, A domain $R$
is $L_1$-hereditary if and only if $R$ is an APD.

$(3)$\ \ Let $R$ be a Noetherian domain. Then $R$ is an APD if and only if $\dim(R)\lst 1$ by {\rm \cite[Theorem 90]{Ka74}.}

$(4)$\ \ Let $R=R_1\times\cdots\times R_n$. Then $R$ is $L_n$-hereditary if and only if each $R_i$ is $L_n$-hereditary.
\end{exa}

\begin{thm}\label{t609}
Let $n\geq 1$. Then the following statements are equivalent for a ring $R$:

$(1)$\ \ $R$ is $L_{n}$-hereditary.

$(2)$\ \ $\pd_{R}M\leq 1$ for any $M \in \mathcal{F}_{n}$.

$(3)$\ \ $\lid_{R}M\leq 1$ for any $M\in \mathcal{F}_{n}$.

$(4)$\ \ Every quotient module of an injective module is $L_{n}$-injective.

$(5)$\ \ $E(M)/M$ is $L_{n}$-injective for any $R$-module $M$, where $E(M)$ is an injective envelope of $M$.

$(6)$\ \ $E(M)/M$ is $L_{n}$-injective for any $R$-module $M$, where $E(M)$ is an $L_{n}$-injective envelope of $M$.

$(7)$\ \ Every submodule $N\in \mathcal{F}_{n-1}$ of a projective $R$-module $P$ is projective.
\end{thm}

\no{\bf Proof.} $(1)\Leftrightarrow (2) \Leftrightarrow (3)$. By Theorem \ref{t504}.

(1)$\Ra$(4). It is trivial.

(4)$\Ra$(1). Let $0\ra N\ra W_0\ra W_1\ra 0$ be exact, where $W_0$ is $L_n$-injective. Let $E$ is the injective hull of $W_0$ and write $W=E/N$.
Then the following is a commutative diagram with exact rows
$$\xymatrix@R=16pt@C=25pt{
0\ar[r]&N\ar[r]\ar@{=}[d]&W_0\ar[r]\ar[d]&W_1\ar[r]\ar[d]&0\\
0\ar[r]&N\ar[r]&E\ar[r]&W\ar[r]&0}$$
Therefore, $0\ra W_0\ra E\bigoplus W_1\ra W\ra 0$ is exact. By hypothesis, $W$ is $L_n$-injective.
By Proposition \ref{p205} and Example \ref{e203}, $W_1$ is $L_n$-injective.

(1)$\Ra$(6)$\Ra$(5). It is trivial.

(5)$\Ra$(4). Let $0\ra K\ra E\ra C\ra 0$ be exact, where $E$ is injective.  Set $E(K)\subseteq E$ is the injective envelope of $K$. Then there exist an
$R$-module $E_{0}$ such that $E=E(K)\bigoplus E_{0}$. So $C\cong E/K\cong (E(K)\bigoplus E_{0})/K\cong (E(K)/K)\bigoplus E_{0}$. Hence $C$
is $L_{n}$-injective since $E(K)/K$ is $L_{n}$-injective by hypothesis.

(4)$\Ra$(7). Let $P$ be a projective $R$-module and $N\in \mathcal{F}_{n-1}$ be a submodule of $P$ and $X$ be any $R$-module. Then there exist exact
sequences $0\ra N\ra P\ra P/N\ra 0$ with $P/N\in \mathcal{F}_{n}$ and $0\ra X\ra E\ra C\ra 0$ with $E$ injective. By hypothesis, $C$ is $L_{n}$-injective.
For any $\alpha \in \Hom_{R}(N, C)$, consider
the following commutative diagram with exact rows
$$\xymatrix@R=16pt@C=25pt{
0\ar[r]&N\ar[r]^{f}\ar[drr]^{\alpha}\ar[dr]_{\theta}&P\ar[r]^{g}\ar[dr]_{\beta}\ar[d]^{\gamma}&P/N\ar[r]&0\\
0\ar[r]&X\ar[r]^{f'}&E\ar[r]^{g'}&C\ar[r]&0}$$
Then there exist $\beta\in \Hom_{R}(P, C)$ such that $\alpha=\beta f$ by Proposition \ref{p204}.
So there exist $\gamma\in \Hom_{R}(P, E)$
such that $\gamma=g'\beta$ since $P$ is projective. There exist $\theta=\gamma f
\in \Hom_{R}(N, E)$ such that $\alpha=g'\theta$.
Then $ \Hom_{R}(N, E)\ra \Hom_{R}(N, C)\ra 0$ is exact. Hence $N$ is projective since $\Ext_{R}^{1}(N,X)=0$.

(7)$\Ra$(1). Let $W$ be an $L_{n}$-injective $R$-module and $N$ be its submodule and $A \in \mathcal{F}_{n}$. Then there exist an
sequence $0\ra K\ra P\ra A\ra 0$ with $P$ projective and $K\in \mathcal{F}_{n-1}$. By hypothesis, $K$ is projective. For any $\alpha \in \Hom_{R}(K, W/N)$, consider
the following commutative diagram with exact rows
$$\xymatrix@R=16pt@C=25pt{
0\ar[r]&K\ar[r]^{f}\ar[drr]^{\alpha}\ar[dr]_{\beta}&P\ar[r]^{g}\ar[d]^{\gamma}&A\ar[r]&0\\
0\ar[r]&N\ar[r]^{f'}&W\ar[r]^{g'}&W/N\ar[r]&0}$$

Then there exist $\beta\in \Hom_{R}(K, W)$ such that $\alpha=g'\beta$. By hypothesis, $W$ is $L_{n}$-injective and $A\in \mathcal{F}_{n}$, then exist
$\gamma \in \Hom_{R}(P, W)$ such that $\beta=\gamma f$ by Proposition \ref{p204}. Then $ \Hom_{R}(P, W/N)\ra \Hom_{R}(K, W/N)\ra 0$ is exact. Hence $W/N$
is $L_{n}$-injective since $\Ext_{R}^{1}(A,W/N)=0$.
\hfill$\Box$

\begin{cor}\label{c610}
Let $n\geq 1$. Let $R$ be a commutative ring and let $S$ be a multiplicatively closed subset of $R$. If $R$ is an $L_n$-hereditary ring,
then so is $R_S$.
\end{cor}

\no{\bf Proof.} Let $M$ be an $R_S$-module with $\fd_{R_S}M\leq n$. Then $\fd_RM\leq n$. By Theorem \ref{t609},
$\pd_RM\leq 1$. Certainly, $\pd_{R_S}M\leq 1$. Consequently, $R_S$ is $L_n$-hereditary.\hfill$\Box$

\begin{cor}\label{t611}
Then the following statements are equivalent for a ring $R$:

$(1)$\ \ $l.L_1\dim(R)\leq 1$.

$(2)$\ \ Every quotient module of an injective module is $L_1$-injective.

$(3)$\ \ Every flat submodule of a projective module is projective.

$(4)$\ \ If $\fd_RM\leq 1$, then $\pd_RM\leq 1$.
\end{cor}

\begin{thm}\label{t612}
Let $n>1$. Then the following statements are equivalent for a ring $R$:

$(1)$\ \ $l.\ld(R)\leq 1$.

$(2)$\ \ $l.L_2\dim(R)\leq 1$.

$(3)$\ \ $l.\FPD(R)\leq 1$.

$(4)$\ \ $l.L_1\dim(R)\leq 1$ and $l.\FFD(R)\leq 1$.
\end{thm}

\no{\bf Proof.} $(1)\Ra (2)$. It is trivial.

(2)$\Ra$(3). Let $M$ be an $R$-module with $k:=\pd_RM<\infty$. If $k>1$, then it is certain that there is an $R$-module
$N$ with $\pd_RN=2$. By hypothesis $\pd_RN\leq 1$, a contradiction. So we obtain $k\leq 1$, whence $l.\FPD(R)\leq 1$.

(3)$\Ra$(1). If $M\in {\mathcal F}_n$, then $\pd_RM<\infty$ by \cite[Proposition 6]{Je70}. Thus $\pd_RM\leq 1$ by hypothesis.
Hence $l.\ld(R)\leq 1$.

$(1)\Ra (4)$. Certainly, $l.L_1\dim(R)\leq l.\ld(R)\leq 1$. For any $M\in {\mathcal F}_n$, we have
$\fd_RM\leq \pd_RM\leq 1$ by Theorem \ref{t609}. Hence ${\mathcal F}_1={\mathcal F}_n$. By Theorem \ref{t211}, ${\rm FFD}(R)\leq 1$.

(4)$\Ra$(1). Since $l.\FFD(R)\leq 1$, ${\mathcal L}_1={\mathcal L}_n$ by using Theorem \ref{t211} again. Hence we have
$l.\ld(R)=l.L_1\dim(R)\leq 1$.\hfill$\Box$

By \cite[Proposition 6]{Je70} and \cite[Corollary 6.4]{FL09}, we have the following corollary:

\begin{cor}\label{cor612a}
A domain $R$ is an APD if and only if $\FPD(R)\leq 1$.
\end{cor}

\begin{thm}\label{t613}
Let $n\geq 1$. Then the following statements are equivalent for a ring $R$:

$(1)$\ \ $R$ is left hereditary.

$(2)$\ \ $R$ is $L_n$-hereditary and $\wgl(R)\leq 1$.

$(3)$\ \ $R$ is $L_1$-hereditary and $\wgl(R)\leq 1$.

$(4)$\ \ $R$ is $L_n$-hereditary and $\wgl(R)<\infty$.

$(5)$\ \ $R$ is $L_n$-hereditary and $\wgl(R)\leq n$.
\end{thm}

\no{\bf Proof.} $(1)\Rightarrow (2)\Rightarrow (3)$ and $(2)\Ra(4)$ are trival.

(3)$\Ra$(1). By Theorem \ref{t601}, every $L_1$-injective module is injective. Hence $R$ is hereditary.

(5)$\Ra$(1). Be similar to (3)$\Ra$(1).

$(4)\Ra(5)$. It is clear for the case $n=1$ by the argument above. Now we let $n>1$. If $k:=\wgl(R)>n$, then there is an
$M\in {\mathcal F}_k$. Let $B$ be the $(k-n)$-syzygy of $M$, then $\fd_RB\lst n$.
By Theorem \ref{t609}, $\fd_RB\leq \pd_RB\leq 1$. Hence $k=\fd_RM\leq k-n+1$. Thus $n\leq 1$, a contradiction. So $k\leq n$,
which implies $\wgl(R)\leq n$. \hfill$\Box$

\begin{thm}\label{t614}
Let $R$ be a Noetherian domain with $dim(R)\leq 1$. Then $R$ is $L_2$-hereditary.
\end{thm}

\no{\bf Proof.} By \cite{Gr71}, $\FPD(R)=\dim(R)\leq 1$. Hence $R$ is $L_2$-hereditary by Theorem \ref{t612}.\hfill$\Box$

\vskip5mm
\section{Examples}

\vskip3mm
Let $M$ be an $R$-module. We say that $M$ is torsion-free if $ux=0$ implies $x=0$, where $u$ is a non-zero-divisor of $R$ and $x\in M$.
It is well-known that a flat $R$-module is certainly torsion-free.

\begin{lem}\label{l701}
Let $R$ be a commutative ring and let $u\in R$ be neither a zero divisor nor a unit. Set $\ol{R}=R/Ru$.
If $A$ be a nonzero $\ol{R}$-module with $\pd_{\ol{R}}A<\infty$, then $\pd_RA=\pd_{\ol{R}}A+1$.
\end{lem}

\no{\bf Proof.}\ \ Let $k=\pd_{\ol{R}}A$. Because of $\pd_R{\ol{R}}=1$, we get $\pd_RA\lst \pd_{\ol{R}}A+1=k+1$ by chang Theorem of rings.
If $k=0$, then we have $\pd_RA=1$ since $A$ is not torsion-free. Hence the assertion holds for $k=0$.

Let $k>0$. By \cite[Exercise 9.6]{Ro79} there is a free $R$-module $F$ such that $\Ext_{\ol{R}}^k(A,F/aF)\not=0$. By Rees Theorem (see
\cite[Theorem 9.37]{Ro79}, we have $\Ext_R^{k+1}(A,F)\not=0$. Hence $\pd_RA\gst k+1$. Thus we get $\pd_RA=k+1$.\hfill$\Box$

\begin{lem}\label{l702}
Let $R$ be a domain and let $J$ be an ideal of $R$ generated by a regular sequence $u_1,\cdots,u_n$. Then
we have{\rm :}

$(1)$\ \ If $M$ is a nonzero $R/J$-module with $\pd_{R/J}M<\infty$, then $\pd_RM=\pd_{R/J}M+n$.

$(2)$\ \ $\pd_RR/J=n$.

$(3)$\ \ For all $k<n$, $\Ext_R^k(R/J,R)=0$.

$(4)$\ \ For all $R/J$-modules $M$ and any $k<n$, $\Ext_R^k(M,R)=0$.

$(5)$\ \ Let $R$ be coherent and set $T=\Ext_R^n(R/J,R)$. Then $\Ext_R^n(T,R)\cong R/J$. Therefore, $T\not=0$.
By the way, $\pd_RT=n$.

$(6)$\ \ If $C$ is a $(n-1)$-cosyzygy of $R$, then $\Ext_R^1(R/J,L_{n-1}(C))\not=0$.
\end{lem}

\no{\bf Proof.}\ \ (1)\ \ Set $R_1=R/(u_1)$. If $n=1$, the assertion holds by Lemma \ref{l701}. Now we assume $n>1$. Then
$\ol{a}_2,\cdots,\ol{a}_n$ is a regular sequence in $R_1$. Thus we may assume by induction that $\pd_{R_1}M=\pd_{R/J}M+(n-1)$.
Using Lemma \ref{l701} again we get $\pd_RM=\pd_{R_1}M+1=\pd_{R/J}M+n$.

(2)\ \ This is direct from (1) by taking $M=R/J$.

(3)\ \ Since $R/J$ is a torsion module, $\Ext^0_R(R/J,R)=\Hom_R(R/J,R)=0$ for $n\geq 1$. Let $n>1$. Then
$\Ext_R^1(R/J,R)=\Hom_{R_1}(R/J,R_1)=0$ by Ress Theorem. Therefore, the assertions for $n=1$ and $n=2$ hold. Let $n>2$ and assume
by induction $\Ext_{R_1}^k(R/J,R_1)=0$ for $k<n-1$. By Rees Theorem we obtain $\Ext_R^{k+1}(R/J,R)=\Ext_{R_1}^k(R/J,R_1)=0$.

(4)\ \ Let $0\ra A\ra F\ra M\ra 0$ be an exact sequence, where $F$ is a free $R/J$-module, that is, $F=\bigoplus R/J$. Thus we have
$\Ext_R^k(F,R)=\prod\Ext_R^k(R/J,R)=0$ by (1). By induction we assume $\Ext_R^{k-1}(A,R)=0$. From the
exact sequence $\Ext_R^{k-1}(A,R)\ra \Ext_R^k(M,R)=\Ext_R^k(F,R)=0$ we have $\Ext_R^k(M,R)=0$.

(5)\ \ Let
$$0\ra P_n\ra P_{n-1}\ra\cdots \ra P_1\ra R\ra R/J\ra 0\eqno{(7.1)}$$
be a projective resolution of $R/J$ in which each $P_i$ is finitely generated. By taking the dual and using the facts $\Ext_R^k(R/J,R)=0$
for $k<n$ we obtain the following resolution
$$0\ra R^*\ra P_1^*\ra\cdots\ra P_{n-1}^*\ra P_n^*\ra T\ra 0.\eqno{(7.2)}$$
Note that $T$ is an $R/J$-module. From (1) we have $\Ext_R^k(T,R)=0$ for $k<n$. By double dual we have the following exact sequence
$$0\ra P_n\ra P_{n-1}\ra \cdots\ra P_1\ra R\ra \Ext_R^n(T,R)\ra 0.$$
It follows that $\Ext_R^n(T,R)\cong R/J$.

From the resolution (7.2) and the fact $\Ext_R^n(T,R)\cong R/J$ we obtain $\pd_RT=n$.

(6)\ \ Note that $\Ext_R^1(R/J,C)\cong \Ext_R^n(R/J,R)\not=0$. Consider the exact sequence $0\ra C\ra L_{n-1}(C)\ra B\ra 0$ in (3.1). Then
we have the following exact sequence
$$\Ext_R^n(R/J,L_n(C))\ra \Ext_R^n(R/J,C)\ra \Ext_R^{n+1}(R/J,B)=0.$$
Hence $\Ext_R^1(R/J,L_{n-1}(C))\not=0$.\hfill$\Box$

\begin{exa}\label{e703}
Now we exhibit a ring $R$ in which ${\mathcal L}_0\supset{\mathcal L}_1\supset {\mathcal L}_2\supset \cdots \supset {\mathcal L}_n\supset\cdots$. To do this, we take
$F$ be a field and set $R=F[x_1,x_2,\cdots,x_n,\cdots]$, where $x_1,x_2,\cdots,x_n,\cdots$ are indeterminates over $F$. Then $R$ is a coherent
domain. For any $n\geq 1$, set $J=(x_1,x_2,\cdots,x_n)$. Then $\Ext_R^n(R/J,R)\not=0$ by Lemma {\rm \ref{l702}}. Let $C_{n-1}$ be a $(n-1)$-cosyzygy
of $R$. By using Lemma {\rm \ref{l702}} again, $L_{n-1}(C_{n-1})$ is an $L_{n-1}$-injective module
but not $L_n$-injective. Thus we are done.
\end{exa}

\begin{exa}\label{e704}
There is a ring $R$ that is $L_1$-hereditary but not almost perfect. In fect, let $D$ be an APD but not a field. Then $D$ is not perfect.
Thus $R=D\times D$ is an $L_1$-hereditary ring and $I=(D,0)$ is a nonzero ideal of $R$. Then $D\cong R/I$ is a proper epic image of $R$
but not perfect. Hence $R$ is not almost perfect.
\end{exa}

\no\textbf{Acknowledgements}

The second author was supported by the National Natural Science
Foundation of China (Grant No. 11171240) and by the Research
Foundation for Doctor Programme (Grant No. 20125134110001). The
fifth author was supported by the National Natural Science
Foundation of China (Grant No. 11401493).

\end{document}